\documentclass
{amsart}
\usepackage{amssymb,amsmath,amsthm,graphicx,pb-diagram,lamsarrow,pb-lams,psfrag
}
\def\lbl#1{\label{#1}\printname{#1}}
\theoremstyle{plain}

\newtheorem{theorem}{Theorem}
\newtheorem{proposition}{Proposition}[section]
\newtheorem{lemma}[proposition]{ Lemma}

\theoremstyle{definition}
\newtheorem{definition}[proposition]{ Definition}

\newtheorem{example }[proposition]{Example}

\newtheorem{remark}[proposition]{Remark}

\newcommand{\la}{\langle}
\newcommand{\ra}{\rangle}

\newcommand{\noi}{\noindent}

\renewcommand{\hom}{\operatorname{Hom}}
\renewcommand{\ker}{\operatorname{ker }}

\newcommand{\im}{ \operatorname{Im}}
\newcommand{\con}{\equiv}

\renewcommand{\sp}{\operatorname{Sp}}

\newcommand{\psdraw}[3]{\centerline{\includegraphics[height=#1,width=#2]{#3}}}
\def\printname#1{\if\draft y\smash{\makebox[0pt]{\hspace{-0.5in}\raisebox{8pt}{\tt\tiny #1}}}\fi
}

\def\cd {\cdot}
\def\ld{\ldots}
\def\LL{\Lambda}

\def\beq{ \begin{eqnarray*}}
\def\enq{\end{eqnarray*}}
\def\en{\end{array}\right]}
\def\a{\alpha} 
\def\b{\beta}
\def\d{\delta}

\def\g{\gamma}

\def\md{\mathcal D}
\def\Z{\mathbb Z}

\def\D{\mathsf D}
\def\H{\mathcal H}
\def\C{\mathcal C}

\def\A{\mathcal A}
\def\B{\mathcal B}

\def\L{{\mathsf L}}

\def\G{\mathcal G}
\def\Q{\mathbb Q}

\def\M{\mathcal M}

\def\bd{\partial}
\def\sg{\Sigma_{g,1}}

\def\iso{\cong}

\def\sub{\subseteq}

\def\sups{\supseteq}
\def\con{\equiv}

\def\pgf{\mathcal P_g^{\text{fr}}}

\def\mg{\mathcal M_{g,1}}

\def\Hg{H_g}

\def\ssg{\Sigma_g}

\def\T{\mathcal T}

\def\tg{T_g}

\def\bg{\mathcal B_g}

\def\ti{\tilde}

\newcommand{\mgg}[1]{\mathcal M_{#1 ,1}}
\newcommand{\fk}{\mathcal F_k (\mathcal M_{g,1})}
\newcommand{\lk}{\mathcal F^L_k (\mathcal M_{g,1})}
\newcommand{\gk}{\mathcal G_k (\mathcal M_{g,1})}
\newcommand{\ffk}[2]{\mathcal F_{#1} (\mathcal M_{#2,1})}
\newcommand{\ggk}[2]{\mathcal G_{#1} (\mathcal M_{#2,1})}
\newcommand{\flk}[2]{\mathcal F^L_{#1}(\mathcal M_{#2,1})}
\newcommand{\glk}[2]{\mathcal G^L_{#1} (\mathcal M_{#2,1})}

\def\draft{n}
\begin{document}
\title
[Lagrangian filtration and homology spheres]{\bf The
Lagrangian filtration of the mapping class group and finite-type invariants of homology spheres}
\author{Jerome Levine}
\address{Department of Mathematics, 
         Brandeis University \newline\indent
         Waltham, MA 02454-9110, USA}        
\email{levine@brandeis.edu  }
\urladdr{http://people.brandeis.edu/$\!\sim$levine/} 
\date{\today}
\thanks{Partially supported by an NSF grant
           and by an Israel-US BSF grant.}
\keywords{Homology spheres, mapping class group }
\subjclass[2000]{Primary (57N10); Secondary (57M25)}.
\begin{abstract}
In a recent paper we defined a new filtration of the mapping class group---the {\em Lagrangian} filtration. We here determine the successive quotients of this filtration, up to finite index. As an application we show that, for any additive invariant of finite-type (e.g. the Casson invariant), and any level of the Lagrangian filtration, there is a homology 3-sphere which has a Heegaard decomposition whose gluing diffeomorphism lies  at that level, on which this invariant is non-zero. In a final section we examine the relationship between the Johnson and Lagrangian filtrations.
\end{abstract}
\maketitle
\section{Introduction}
Dennis Johnson in \cite{J} defined a filtration of the mapping class group of a closed orientable surface $\Sigma$ based upon the induced action of its elements on the lower central series of the fundamental group of $\Sigma$. He also showed how to imbed the graded Lie algebra formed by the successive quotients of this filtration into a graded Lie algebra $H\otimes \L (H)$, where $H=H_1 (\Sigma )$ and $\L (H)$ is the free Lie algebra on $H$. It has been a long-standing problem to determine the image of this imbedding.

More recently we defined another filtration of the mapping class group, the Lagrangian filtration \cite{L}, and defined an analogous imbedding of the graded group formed by the successive quotients into $H'\otimes\L (H')$, where $H'=H_1 (T)$ and $T$ is the handlebody bounded by $\Sigma$. A major difference between the two filtrations is that the Johnson filtration has trivial residue (i.e. intersection of the sequence of filtrations), and so the totality of successive quotients approximates the full mapping class group, while the residue of the Lagrangian filtration is the subgroup consisting of diffeomorphisms which extend to diffeomorphisms of $T$. Therefore the successive quotients of the Lagrangian filtration do not approximate the the full mapping class group, but they do serve as well for the purpose of studying Heegaard decompositions of 3-manifolds.

The main result of this paper is  that the image  in $H'\otimes\L (H')$ of the successive quotients of the Lagrangian filtration is a subgroup of finite index at each level, if the genus of $\Sigma$ is large enough. The analogous assertion is definitely known to be false for the Johnson filtration (although it is true, and actually onto $H\otimes\L (H)$, see  \cite{L}, for the Johnson filtration of the group of {\em homology cylinders} over $\Sigma$).

As an application we consider the following question. Given a homology 3-sphere $M$, how far down in a given filtration of the Torelli group can one find a diffeomorphism such that $M$ has a Heegaard decomposition with that gluing diffeomorphism? For the Johnson filtration it is now known that this can be done down to the third level (see \cite{Mo}, \cite{P}). We consider this question for the Lagrangian filtration and conjecture that it can be done arbitrarily far down. As evidence for this conjecture we show that, for any additive rational invariant of finite type $\lambda$, for example the Casson invariant, there exist $h$ arbitrarily far down in the Lagrangian filtration such that, for the associated homology sphere $M_h$, $\lambda (M_h)\not= 0$. (It follows from S. Morita's fundamental work  \cite {Mo} and the deep results of R. Hain \cite{Ha},  that, for the Casson invariant $\lambda$, there exist $h$ arbitrarily far down, even in the Johnson filtration, such that $\lambda (M_h)\not= 0$ (see \cite[Theorem 14.10]{Ha}), but it is unknown whether this is true for other finite-type invariants.)

In the final section we examine the question of whether the Lagrangian filtration might be just the product of the Johnson filtration and the residue of the Lagrangian filtration. We present a recursive approach to this question and use this approach, together with results of Johnson and Morita on the first Johnson homomorphism, to show that it is true in low degrees.

\section{The Johnson filtration}
Let's begin with a quick outline of the basic facts about the Johnson filtration of the mapping class group. For more details and additional references see \cite{J} and \cite{M}. 

We will restrict our attention to the mapping class group $\mg$, defined to be the
group of isotopy classes,  rel boundary, of  diffeomorphisms of the oriented
surface $\sg$ of genus $g$ with one boundary component, which are the identity on
the boundary. Throughout this paper diffeomorphism will mean orientation-preserving diffeomorphism. By a classical theorem of Dehn-Nielsen-Baer, $\mg$ can be
identified with the group of automorphisms of $F=F(x_i ,y_i )=\pi_1 (\sg )$, the
free group with generators $x_1 ,\ld ,x_g ,y_1 ,\ld ,y_g$ corresponding to the
meridian and longitudinal curves of $\sg$, which fix the element $\d =[x_1 ,y_1
]\cdots [x_g ,y_g ]$. Johnson, in \cite{J}  proposed the following filtration of $\mg$. Let
$\fk$ be the subgroup of $\mg$ consisting of all $h$ which induce the identity
automorphism of $F/F_{k+1}$. $F_q$ denotes the lower central series subgroup of
$F$ generated by all commutators of weight $q$ or more. For example $\ffk{1}{g}$ is the same as the
Torelli group $\T_g$. Since $F$ is residually nilpotent $\cap_k \fk=\{1\}$. If $h\in\fk$ then the equality $h_*
(g)=g\bar h (g)$, for any $g\in F$, defines a function $\bar h:F\to F_{k+1}$,
which then induces a homomorphism $F/F_2\to F_{k+1}/F_{k+2}$. Set $\Hg =F/F_2
=H_1 (\sg )$. Recall that the lower central series quotients $\{ F_q
/F_{q+1}\}$ form a graded Lie algebra, where the Lie bracket is induced by the
commutator in $F$. In fact this Lie algebra is isomorphic to the free Lie algebra $\L(\Hg )$.
The assignment $h\to\bar h$ defines a function $J_k :\fk\to\hom (\Hg,
\L_{k+1}(\Hg ))$, which is, in fact, a homomorphism. It is obvious that $\ker J_k
=\ffk{k+1}{g}$. The identification of the image of $J_k$ remains one of the  major
problems in the study of the algebraic structure of the mapping class group. 

The commutator in $\mg$  induces a Lie algebra structure on the graded
abelian group $\{\gk =\fk /\ffk{k+1}{g}\}$. There is also a Lie algebra structure on
the graded abelian group $\{\hom (\Hg ,\L_{k+1}(\Hg ))\}$, which arises from the observation that any homomorphism $\Hg\to \L_{k+1}(\Hg)$ can be uniquely extended to a
derivation of $\L(\Hg)$ which is of degree $k$, and thus we can identify  $\hom (\Hg, \L_{k+1}(\Hg))$ with $\md_k
(\L(\Hg))$, the degree $k$ component of  $\md (\L (\Hg ))$, where $\md(L)$ denotes the graded Lie algebra of derivations of a Lie
algebra $L$. Furthermore the collection $\{ J_k\}$ defines a Lie algebra homomorphism
$\G (\mg )\to\md (\L(\Hg))$ and so the image is a Lie subalgebra of $\md (\L(\Hg))$.

In \cite{M} S. Morita pointed out that $\im J$ is contained in a natural Lie
subalgebra of $\md (\L(\Hg))$, denoted $\D (\Hg)$, which is defined as follows. First of all the
symplectic form on $\Hg$ induces a canonical isomorphism $\hom (\Hg, \L_{k+1}(\Hg))\iso \Hg\otimes
\L_{k+1}(\Hg)$. If $\{ x_i ,y_i\}$ is a symplectic basis of $\Hg$ then the
isomorphism can be explicitly defined by $\phi\to\sum_i (x_i\otimes\phi (y_i )-y_i\otimes\phi
(x_i ))$ but it is independent of the particular choice of basis and depends only upon the symplectic form. Then $\D_k (\Hg)$ is defined to be the kernel of the bracket map $\b_k
:\Hg\otimes \L_{k+1}(\Hg)\to \L_{k+2}(\Hg)$. It is not hard to see that $\{\b_k\}$
coincides with the map $\b :\md (\L(\Hg))\to \L(\Hg)$ defined by $\b (d)=d(\sum_i [x_i
,y_i ])$. $\b$ is not a map of Lie algebras but it is easy to check that $\ker\b$
is a Lie subalgebra of $\md (\L(\Hg))$ (but not an ideal).

In \cite{Jo} Johnson pointed out that $\D_1 (\Hg)$ could be identified with the
exterior power $\LL^3 \Hg$ and showed that $\im J_1 =\D_1 (\Hg)$. In \cite{Mo} Morita
showed that $\im J_2$ is a subgroup of $\D_2 (\Hg)$ of index some power of $2$.
On the other hand Morita showed in \cite{M} that $\im J_k$ is a subgroup of
infinite index, if $k$ is odd $>1$, by constructing a homomorphism $\hom (\Hg, \L_{k+1}(\Hg))\to P^k (H_g )$, the {\em Morita trace}, where $P^k (H_g )$ is the $k$-th symmetric power of $H_g$. This trace function is
$0$ on $\im J_k$ for $k>1$ but, for $k$ odd, is rationally onto, even when restricted to $\D_k (H_g )$, (but $0$ for $k$ even). Finally, for $k$ even $>2$
work of Nakamura \cite{N} shows that $\im J_k$ is generally a subgroup of infinite index.

\section{The Lagrangian filtration}

In \cite{L} I introduced a different filtration on $\mg$ which I will recall
here. Note that there is a change of notation from \cite{L} and I will discuss only one of the
two variations which were treated there (in fact the other one seems to
have some problems). 

Let $\tg$ be a handlebody bounded by $\ssg$, the closed
oriented surface of genus $g$, and let $i:\sg\to\tg$ be the inclusion. Choose
a basis $\{ x_i ,y_i\}$ of $\pi_1 (\sg )$ representing a system of meridian and longitudinal curves so that the $\{ x_i\}$ are null-homotopic
in $\tg$ and therefore  $\{i_* (y_i )\}$ are a basis for $F'=\pi_1 (\tg )$. 
\begin{definition}\lbl{def.lag}
For  $k\ge 1$, define $\lk$ to be the set of all $h\in\mg$ satisfying
\begin{enumerate}
\item $i_* h_* (x_i )\in
F'_{k+1}$ 
\item $h_* (x_ i)\con x_i \mod F_2$ for all $i$
\end{enumerate}
\end{definition}
 We note the following facts, which were proved in \cite{L}
\begin{enumerate}
\item $\lk$ is a subgroup.
\item $\flk{1}{g}$ contains the Torelli group.
\item $\flk{2}{g}$ is the subgroup of $\mg$ generated by Dehn twists along simple closed curves which bound in $\tg$.
\item $\flk{\infty}{g}=\cap_k \lk$ coincides with $\bg\cap\flk{1}{g}$, where $\bg$ is the subgroup of $\mg$ consisting of those diffeomorphisms  which extend to diffeomorphisms of $\tg$.
\item $\flk{k+1}{g}$ is a normal subgroup of $\lk$, the kernel of a homomorphism $J_k^L :\lk\to\D_k (\Hg')$, where $\Hg' =H_1(\tg )$.
\end{enumerate}
To elaborate on item (5), we first define 
$$J_k^L :\lk\to\hom (L,\L_{k+1}(\Hg'))\iso \Hg'\otimes \L_{k+1}(\Hg')$$
 where $L=\ker\{ \Hg\to \Hg'\}$ and the isomorphism is induced by the symplectic form on $\Hg$. If $h\in\lk$ and $l\in L$, then we choose $\lambda\in\ker i_*\sub F$ representing $l$ and define $J_k^L (h)\cd l$ to be the reduction of $i_* h_* (\lambda )\in F'_{k+1}/F'_{k+2}\iso \L_{k+1}(\Hg')$. It is not hard to see that this reduction of $i_* h_* (\lambda )$  depends only on $l$. Then $J_k^L$ is a homomorphism---this requires property (2) in Definition \ref{def.lag}. Clearly $\ker J_k^L =\flk{k+1}{g}$. Thus we obtain an imbedding $\glk{k}{g}\sub \Hg'\otimes \L_{k+1}(\Hg')$ and it is shown in \cite{L} that the image of $J_k^L$ is contained in $\D_k (\Hg')$.

\subsection{Relation to homology spheres}\lbl{sec.hs}

The Lagrangian filtration says less about the structure of $\mg$ than the Johnson filtration, since $\flk{\infty}{g}$ is non-empty, but is more relevant to the classification of homology $3$-spheres via the Heegaard decomposition, since $\flk{\infty}{g}$ consists entirely of diffeomorphisms associated to a Heegaard decomposition of $S^3$. 

Suppose we consider an explicit correspondence which associates to any $h\in\mg$ the oriented $3$-manifold $M_h =\tg\cup_{r_g h}\tg =\tg\amalg\tg$, where every $x\in\bd\tg=\ssg$ (the left hand copy) is identified with $r_g h(x)\in\bd\tg$ (the right-hand copy) and oriented consistent with the left-hand copy of $\tg$. $h$ is extended over $\ssg$  by the identity on the attached disk and $r_g$ is an involution of $\ssg$ which exchanges the meridian and longitude curves.  We can specify $r_g$ by its action on $\pi_1 (\sg )$:
$$r_g (x_i )=x_i y_i x_i^{-1} ,\quad r_g (y_i )=x_i^{-1}$$
It is not hard to see that $M_{\text id}=S^3$ and that $M_h$ is a homology sphere if $h_* (L)=L$ and, therefore, for any $h\in\flk{1}{g}$. It is a standard fact that the  diffeomorphism class of $M_h$ depends only on the class of $h$ in $\mg/\bg$, the space of left cosets $\{ h\bg\}$.  

Now define a left action of $\bg$ on $\mg/\bg$ by the formula:
$$ b\cd [h]=[r_g^{-1} br_g h]$$
Note that $r_g^{-1} br_g$ is a diffeomorphism of $\ssg$ which extends to a diffeomorphism of the ``dual'' handlebody $\bar\tg$ bounded by $\ssg$, chosen so that the kernel of the inclusion $\pi_1 (\sg )\to\pi_1 (\bar\tg )$ is normally generated by $y_1 ,\cdots ,y_g$. Alternatively $\bar\tg$ can be identified with $\overline{S^3 -\tg}$, where $\tg$ is imbedded in $S^3$ in the standard ``unknotted'' way. In fact the subgroup of $\mg$ consisting of all such diffeomorphisms is exactly $r_g^{-1} \bg r_g$. It is readily apparent now that the  diffeomorphism class of $M_h$ depends only on the orbit of $h$ under this action of $\bg$ on $\mg/\bg$.

These observations hold stably. Consider the canonical suspension $\mg\sub\mgg{g+1}$, defined by imbedding $\sg\sub\Sigma_{g+1,1}$ and extending a diffeomorphism of $\sg$ over $\Sigma_{g+1,1}$ by the identity on $\Sigma_{g+1,1}-\sg$. Then the diffeomorphism class of $M_h$, for $h\in\mg$, is unchanged if we replace $h$ by its suspension in $\mgg{g+1}$. In fact the Reidemeister-Singer theorem can be expressed by saying that $M_h$ and $M_{h'}$ are diffeomorphic if and only if, after sufficient suspensions, $[h']=b\cd [h]$ for some $b\in\bg$. This formulation is also used in  \cite{Mo} and \cite{P}.

With this in mind we are most interested in identifying $\im J_k^L$ stably, i.e. for large $g$.

As an application of our results we will address, in  Section \ref{sec.hss}, the question, raised by Morita in \cite{Mo}, of whether restricting an element $h\in\T_g$ to lie in some suitable subgroup of $\T_g$ imposes some topological restriction on $M_h$. We will show that, for any $\Q$-valued additive finite type invariant $\lambda$ , there are elements $h$ arbitrarily far down in the Lagrangian filtration such that $\lambda$ is  non-zero on $M_h$.

\section{Main results}
In contrast with the  Johnson homomorphisms $\{ J_k\}$, we will show that, for large $g$, every $J_k^L$ is, up to 2-torsion, onto.
\begin{theorem}\lbl{th.main} The image of $J_k^L$ contains $2\D_k (\Hg')$, and so is a subgroup of index a power of $2$, if $g>k$. If, in addition, $k$ is odd, then $J_k^L$ is onto.
\end{theorem} 
\begin{remark}\lbl{rem.Mo}
For $k=2$ we can use Morita's calculation of $J_2$ in \cite{Mo} to show that $J_2^L$ is actually onto. So it is reasonable to ask whether this theorem can be improved to say that $J_k^L$ is onto for all $k<g$. 
\end{remark}
\begin{remark}
In \cite{L} we made use of an imbedding of the framed  pure braid group $\pgf$ into the mapping class group $\mg$, first defined by Oda \cite{O}, to define an imbedding 
$$(\pgf )_{k+1}/(\pgf )_{k+2} \sub\glk{k}{g}$$
This gives a lower bound for the rank of $\glk{k}{g}$ which, as was pointed out in \cite{L}, is generally strictly lower than the rank of $\D_k (\Hg')$. 
\end{remark}

\subsection{Proof of Theorem \ref{th.main}}
We will actually prove a stronger fact. Consider $\Hg'\sub \Hg$ to be the subgroup generated by the $\{ y_i\}$ and so $\D_k (\Hg')\sub\D_k (\Hg)$. We will prove:
\begin{lemma}\lbl{lem.main}
$\im J_k\sups 2\D_k (\Hg')$ if $g>k$. If $k$ is odd then $\im J_k\sups \D_k (\Hg')$.
\end{lemma}
Theorem \ref{th.main} will follow immediately, since $\fk\sub\flk{k}{g}$ and  $J_k^L |\ffk{k}{g}$ is just $J_k$ followed by the projection from $\D_k (\Hg)$ to $\D_k (\Hg')$ (induced by the projection $\Hg\to \Hg'$).

\begin{remark}
The proof will actually show that $J_k ((\T_g )_k )\sups 2\D_k (\Hg')$ and, for odd $k$, $J_k ((\T_g )_k )\sups \D_k (\Hg')$, providing $g>k$. Recall from \cite{M} that $(\T_g )_k\sub\ffk{k}{g}$. 
\end{remark}

Lemma \ref{lem.main} will follow easily from:
\begin{lemma}\lbl{lem.brak}
There exist subgroups $\ti\D_k (\Hg)\sub\D_k (\Hg)$ which satisfy: 
\begin{enumerate}
\item $\ti\D_k (\Hg)\sups 2\D_k (\Hg)$. If $k$ is odd then $\ti\D_k (\Hg)=\D_k (\Hg)$.
\item If $g>k$ then any element $\a\in\ti\D_k (\Hg')=\ti\D_k (\Hg)\cap\D_k (\Hg')$ is a linear combination $\sum_i  [\b_i ,\g_i ]$, where $\b_i\in\D_1 (\Hg)=\ti\D_1 (\Hg)$ and $\g_i\in\ti\D_{k-1}(\Hg')$. 
\end{enumerate}
\end{lemma}
 To prove Lemma \ref{lem.main}  from Lemma \ref{lem.brak} we show  that $\ti\D_k (\Hg')\sub\im J_k$ by induction on $k$.

If $k=1$ this follows immediately from the fact \cite{Jo} that $\im J_1 =\D_1(\Hg)$.

Now suppose $\im J_{k-1}\sups\ti\D_{k-1}(\Hg')$. Let $\a\in\ti D_k (\Hg')$.  By Lemma \ref{lem.brak} we can write $\a=\sum_i  [\b_i ,\g_i ]$ for some $\b_i\in\D_1 (\Hg),\ \g_i\in\ti\D_{k-1}(\Hg')$. Then $\b_i =J_1 (h_i )$ for some $h_i\in\ffk{1}{g}$ and, by induction, $\g_i =J_{k-1}(h'_i )$ for some $h'_i\in\ffk{k-1}{g}$. Now $J$ is a Lie algebra homomorphism $\ggk{}{g}\to\D (\Hg)$ and so
$$\a =\sum_i [\b_i ,\g_i ]=\sum_i [J_1 (\bar h_i ), J_{k-1}(\bar h'_i )]=J_k (\sum_i [\bar h_i ,\bar h'_i ])$$
where $\bar h_i ,\ \bar h'_i$ denote the classes of $h_i ,\ h'_i$ in $\ggk{1}{g},\ \ggk{k-1}{g}$, respectively.

This completes the proof of Lemma \ref{lem.main}.

\begin{proof}[Proof of Lemma \ref{lem.brak}] Recall the {\em quasi-Lie algebra} $\A^t (\Hg )$ of planar binary trees with leaves decorated by elements of $\Hg$, discussed in \cite{GL} and \cite{Le} (see also \cite{HP} and \cite{P} for a similar tree Lie algebra). A quasi-Lie algebra is the same as a Lie algebra except that the axiom $[\a ,\a ]=0$ is replaced by the slightly weaker anti-symmetry axiom $[\a ,\b ]=-[\b ,\a ]$. Thus any quasi-Lie algebra becomes a Lie algebra after inverting $2$. By definition $\A_k^t (\Hg )$ is generated by planar binary trees with $k+2$ leaves whose decorations are members of the basis $\{ x_1 ,\ldots ,x_g ,y_1,\ldots y_g\}$,  subject to what are known as the anti-symmetry and IHX relations. See \cite{Le}  for more details. There is a quasi-Lie algebra homomorphism $\eta:\A^t (\Hg )\to\D(\Hg)$ which assigns to a decorated tree $T$  with $k+2$ leaves an element $\sum_i a_i\otimes\lambda_i\in \Hg\otimes \L_{k+1}(\Hg)$, where the sum ranges over the leaves of $T$, $a_i$ is the decoration of the chosen leaf and $\lambda_i$
is the element of $\L_{k+1}(\Hg)$ associated to the rooted planar binary tree created from $T$ by removing the decoration $a_i$ from the chosen leaf and making that leaf the root. This element lies in $\D_k (\Hg)$. See \cite{Le} or \cite{P} for the details. The Lie bracket in $\A^t (\Hg)$ is defined as follows. If $T_1$ and $T_2$ are decorated trees, then
$$[T_1 ,T_2 ]=\sum_{i,j} \la a_i , b_j\ra T_1\ast_{ij}T_2 $$
The terms of this sum range over all pairs consisting of a leaf of $T_1$, with decoration $a_i$, and a leaf of $T_2$ with decoration $b_j$. $\la ,\ra$ denotes the symplectic pairing on $\Hg$. Then $T_1\ast_{ij}T_2$ denotes the decorated tree constructed by removing the decoration $a_i$ from the $i$-th leaf of $T_1$ and $b_j$ from the $j$-th leaf of $T_2$ and  welding these undecorated leaves together.
See \cite{GL} or  \cite{P} for more details.

It is well-known that $\eta$ induces an isomorphism $\A^t (\Hg )\otimes\Q\iso\D(\Hg)\otimes\Q$. Let us  define $\ti\D (\Hg)=\im\eta$. Then it is proved in \cite{Le} that $\ti\D (\Hg)\sups 2\D (\Hg)$ and that,  if $k$ is odd, $\ti\D_k (\Hg)=\D_k (\Hg)$. Thus $\ti\D (\Hg)$ satisfies $(1)$ of Lemma \ref{lem.brak}.

We now prove $(2)$. We may assume that $\a$ is represented by a tree $T$ with $k+2$ leaves whose decorations are drawn from $\{ y_i\}$. Choose a pair of leaves with a common trivalent vertex. Denote their two decorations $y_r$ and $y_s$. Now remove the two leaves thereby creating a tree $T'$ with $k+1$ leaves, one of which is the remains of the original trivalent vertex and is so far undecorated. We choose a decoration $y_t$ which, because $g>k$, can be assumed to be different from all the $y_i$ which decorate the other leaves of $T'$. Let $\g$ be the element of $\ti\D_{k-1}$ represented by $T'$. We also create a tree $T''$ with three leaves decorated by $y_r , y_s$ and $x_t$ and let $\b\in\D_1 (\Hg)$ be the element represented by $T''$. It now follows easily from the graphical definition of the bracket described above that $[\b ,\g ]=\pm\a$, since only one pair of decorations from $T'$ and $T''$ have a non-zero value under $\la ,\ra$. The sign depends on the cyclic order of the decorations chosen for $T''$.
See Figure \ref{fig.bra}.

\begin{figure}[ht]
\psdraw{1.5in}{2.5in}{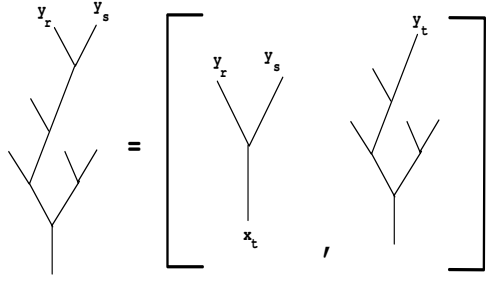}
\bigskip
\caption{Representing a planar binary labeled graph as a bracket.}
\lbl{fig.bra}
\end{figure}

This completes the proof.
\end{proof}
\subsection{Refinement of Theorem \ref{th.main}}
We will need a strengthened form of Theorem \ref{th.main} for our results in Section \ref{sec.hss}. Recall from Section \ref{sec.hs} that $\bar\bg =r_g^{-1}\bg r_g$ is the subgroup of $\mg$ consisting of all diffeomorphisms of $\ssg$ which extend over $\bar\tg$.

\begin{theorem}\lbl{th.main1}
$J_k^L (\bar\bg\cap\flk{k}{g})$ contains $2\D_k (\Hg')$ when $g>k$.

 If, in addition, $k$ is odd, then $J_k^L(\bar\bg\cap\flk{k}{g})=\D_k (\Hg')$.
\end{theorem}
\begin{remark}\lbl{rem.con}
It follows from the results of \cite{P} that $J_2^L(\bar\bg\cap\flk{2}{g})=\D_2 (\Hg')$ for $g>9$. One might conjecture that $J_k^L(\bar\bg\cap\flk{k}{g})=\D_k (\Hg')$ for all $k<g$ or, at least, that $J_k^L(\bar\bg\cap\flk{k}{g})=\im J_k^L$.
\end{remark}
We derive Theorem \ref{th.main1} from the correspondingly strengthened form of Lemma \ref{lem.main}.
\begin{lemma}\lbl{lem.main1}
$J_k (\bar\bg\cap\ffk{k}{g})\sups 2\D_k (\Hg')$ if $g>k$. 

If $k$ is odd then 
$J_k (\bar\bg\cap\ffk{k}{g})\sups\D_k (\Hg')$.
\end{lemma}
Recall that $\Hg'\sub \Hg$ is here taken to be the subgroup generated by the $\{ y_i\}$ and so $\D_k (\Hg')\sub\D_k (\Hg)$.
\begin{remark}\lbl{rem.main1}
The proof will actually show that $J_k (\bar\bg\cap(\T_g )_k )\sups 2\D_k (\Hg')$ and, for odd $k$, $J_k (\bar\bg\cap(\T_g )_k )\sups \D_k (\Hg')$, providing $g>k$. This will be used in the applications in Section \ref{sec.hss}.
\end{remark}

 \begin{proof}[Proof of Lemma \ref{lem.main1}] 
We follow the lines of the proof of Lemma \ref{lem.main}, proceeding  by induction on $k$. For $k=1$ we  use the following result of Morita \cite{Mo}. This is stronger than the assertion of Lemma \ref{lem.main1}, but is needed for the inductive step and lateron, in Section \ref{sec.JvsL}.
 
 \begin{lemma}[\cite{Mo}]\lbl{lem.johnson}
$J_1 (\bar\bg\cap\ffk{1}{g})$ is generated by the elements in $\D_1 (\Hg )\iso\bigwedge^3\Hg$ of the form $y_k\wedge y_l\wedge y_m ,\  x_k\wedge y_l\wedge y_m$ and $x_k\wedge x_l\wedge y_m$ for all $1\le  k,l,m\le g$.
\end{lemma}

For the inductive step in the proof of Lemma \ref{lem.main1} we use  Lemma \ref{lem.brak}. As pointed out in the proof of that Lemma, $\ti\D_k (\Hg')$ is generated by elements represented by trees $T=[T_1 ,T_2 ]$, where $T_1$  represents an element of $\ti\D_{k-1}(\Hg'
)$ and $T_2$ is a tree with three leaves decorated by $y_r ,y_s$ and $x_t$. By induction $T_1$ is realized by an element of $\bar\bg\cap\ffk{k-1}{g}$ and, by Lemma \ref{lem.johnson}, $T_2$ is represented by an element of $\bar\bg\cap\ffk{1}{g}$. Thus $T$ is represented by a product of commutators of elements of $\bar\bg$. 
\end{proof}

This completes the proof of Theorem \ref{th.main1}.
 
\section{Homology spheres}\lbl{sec.hss}
 In considering the relation between the mapping class group $\mg$ and closed $3$-manifolds one can pose the following question. Given any subgroup $G$ of $\mg$ what, if any, topological restriction is imposed on $M_h$ if $h$ is restricted to lie in $G$? Of course the action of $h_*$ on $H_1 (\Sigma_g )$ may impose homological restrictions on $M_h$, so we may refine this question by restricting $G$ to be a subgroup of the Torelli group $\T_g$ and  therefore  $M_h$ to be a homology sphere. We mention three known results
 \begin{enumerate}
\item In \cite{Mo} Morita shows that any homology sphere can be represented as $M_h$ for some $h\in\ffk{2}{g}$ and Pitsch \cite{P} has proved  that this is true for $\ffk{3}{g}$.
\item It is shown in \cite{Mo} that the Casson invariant vanishes on $M_h$ if $h\in (\T_g )_3$ and, more generally in \cite{GL1} and \cite{H} that if $h\in (\T_g )_{k+1}$ then any finite-type invariant of type $k$ vanishes on $M_h$.
\item  In R. Hain's fundamental paper \cite{Ha} on the Malcev Lie algebra of the Torelli group, he shows, using Morita's work on the Casson invariant, that there exists $h$ arbitrarily far down in the Johnson filtration such that the Casson invariant of $M_h$ is non-zero.
\end{enumerate}
 Item (1) prompts the question of whether, for a given $k$,  any homology sphere can be represented as $M_h$ for some $h\in\ffk{k}{g}$. Items (2) and (3) prompt the question of whether, on the contrary, for a given topological invariant $\lambda$, there is some $k$ such that $\lambda (M_h)=0$ for every $h\in\ffk{k}{g}$.
 
 One can pose the same questions for $\flk{k}{g}$. If the conjecture in Remark \ref{rem.con} that $J_k^L(\bar\bg\cap\flk{k}{g})=\im J_k^L$ for all $k$ and sufficiently large $g$ were known to be true, then one can easily  show that, for any given homology sphere $M$ and any positive integer $k$, there exists some $h\in\lk$ (for some $g$) so that $M=M_h$. In the meantime we will  give a partial answer.
 
As a warmup and easy application of our results we first prove the following. As remarked above this result is already known to be true in the much stronger context of the Johnson filtration, as a consequence of deep work of Morita and Hain.
\begin{proposition}\lbl{th.cas}
For any finite $k$  there exists $h\in\flk{k}{g}$ (for some $g$) such that the Casson invariant of $M_h$ is non-zero.
\end{proposition}
\begin{remark}
Proposition \ref{th.cas} is obviously false for $k=\infty$ since $\flk{\infty}{g}\sub\bg$ and if $h\in\bg$ then $M_h =S^3$.
\end{remark}
We give a considerable generalization of this Proposition to the more general class of additive finite-type invariants in the next section.

\begin{proof}
First of all we can, by \cite{Mo}, choose $h\in\ffk{2}{g}$ so that $\lambda (M_h )\not= 0$, where $\lambda$ denotes the Casson invariant. Suppose, by induction, that there exists $h\in\flk{k}{g}\cap\ffk{2}{g}$ so that $\lambda (M_h )\not= 0$. Now $J_k^L (h^2 )=2J_k^L (h)$ and so, by Lemma \ref{lem.main1}, there is some $b\in\bar\bg\cap\ffk{k}{g}$ such that $J_k^L (b)=J_k^L (h^2)$. So $J_k^L (b^{-1} h^2 )=0$ and therefore, since $k\ge 2$, $b^{-1} h^2\in\flk{k+1}{g}\cap\ffk{2}{g}$. But since $M_{b^{-1} h^2}=M_{h^2}$ and, by \cite{Mo}, $h\to\lambda (M_h)$ defines a homomorphism $\ffk{2}{g}\to\Z$ we have
$$\lambda (M_{b^{-1} h^2} )=\lambda (M_{h^2} )=2\lambda (M_h )\not=0$$
This completes the proof.
\end{proof}

\subsection{Finite type invariants of homology spheres}
Before stating and proving our next application we give a brief resum\' e of what we need about rational finite-type invariants. This notion was first defined by Ohtsuki in \cite{Oh}. Let $\H$ denote the $\Q$-vector space with basis the set of diffeomorphism classes of oriented homology 3-spheres. A filtration of $\H$ is defined as follows. Let $L\sub S^3$ be an algebraically split link (i.e. the linking number of any two components of $L$ is $0$) with each component given a $\pm 1$ framing. Consider the element of $\H$ given by
$$\sum_{L'}(-1)^{|L'|}M_{L'}$$
where $L'$ ranges over all sublinks of $L$, $|L'|$ denotes the number of components of $L'$, and $M_{L'}$ is the homology sphere obtained by surgery on $L'$. Then $\H_m$ is the subspace of $\H$ generated by all such elements for all links $L$ with $m$ components. It is shown in \cite{GL1} that $\H_{3m}=\H_{3m+1}=\H_{3m+2}$.

In \cite{GL2} another filtration of $\H$ is defined. Let $\phi$ denote the function which assigns to every $h\in\T_g$ the homology sphere $M_h$. We extend this to a linear function on the group algebra $\Q [\T_g ]\to\H$, also denoted by $\phi$. We then define $\H'_m=\bigcup_g \phi (I_g^m )$, where $I_g$ is the augmentation ideal of $\Q [\T_g ]$. It is shown in \cite{GL2} that $\H'_{2m}=\H'_{2m+1}=\H_{3m}$. Finally the recent theory of claspers of Habiro \cite{H} and Goussarov \cite{Gu}---see also \cite{GGP}---provides another filtration, which turns out to coincide with $\{\H'_m\}$. Habiro defines the notion of $A_k$-surgery equivalence of 3-manifolds and if $M$ and $M'$ are $A_k$-surgery equivalent homology spheres then $M-M'\in\H'_k$. Habiro asserts that $M$ and $M'$ are $A_k$-surgery equivalent if and only if there is a closed oriented surface $\Sigma\sub M$ so that if we cut $M$ open along $\Sigma$ and reglue using some $h\in (\T_g )_k$ ($g=$ genus of $\Sigma$ ), the result  is diffeomorphic to $M'$. In order to relate this fact to ordinary Heegaard decomposition we will need the following lemma. 

Let $i:\Sigma_g\to M$ be an imbedding into the closed oriented 3-manifold $M$, and $h\in\mg$. Define $M_{i,h}$ as follows. Cut $M$ open along $i(\Sigma_g )$ to obtain a manifold $M_0$ with two boundary components, both identified with $\Sigma_g$ by $i$. Then consider $h$ to be a diffeomorphism from the boundary component whose orientation agrees with $M_0$ to the other boundary component and use it to glue the boundary components together to obtain the oriented closed manifold $M_{i,h}$
\begin{lemma}\lbl{lem.hab}
If $i:\Sigma_g\to S^3$ is an imbedding and $h\in\mg$ then $M_{i,h}$ has a Heegaard representation $M_{i,h}=M_{h' }$, where $h'\in\M_{g',1}$ is some conjugate of the suspension of $h$ into $\mgg{g'}$.
\end{lemma}
\begin{proof} We may assume, after some suspensions, that $i(\ssg )$ separates $S^3$ into two components, both of which are handlebodies of genus $g$. To accomplish this first connect sum $i(\ssg )$ with a parallel copy of itself to get $i(\ssg )$ to be separating. Now one component $C'$ of the complement is already a handlebody. To get the other component $C$ to also be a handlebody,  choose a handle decomposition of $C$ based on $\bd C$. We can arrange that there are no $0$-handles and $3$-handles and then adjoin the $1$-handles to $C'$. After this,  the dual handle decomposition of $C$ will have  only  $1$-handles. Note that $C'$ will still be a handlebody.

Now choose an orientation preserving diffeomorphism $t:T_{g'}\to C$, where $C$ now denotes the complementary component whose orientation agrees with $i(\ssg )$ and, similarly, $t':T_{g'}\to C'$ Then we have 
\begin{enumerate}
\item $S^3 =T_{g'}\cup_{{t'}^{-1} t}T_{g'}=M_{r_{g'}^{-1} {t'}^{-1} t}$
\item $M_{i,h}=T_{g'}\cup_{{t'}^{-1} ih'i^{-1} t}T_{g'}=T_{g'}\cup_{{t'}^{-1} tfh'f^{-1}}T_{g'}=M_{r_{g'}^{-1} {t'}^{-1} tfh'f^{-1}}$
\end{enumerate}
where $h'$ is an iterated suspension of $h$ and $f=t^{-1} i$. It follows from (1) and the Reidemeister-Singer theorem, as discussed in Section \ref{sec.hs}, that, after further suspension,  $r_{g'}^{-1} {t'}^{-1} t=r_{g'}^{-1} b_1 r_{g'}b_2$ for some $b_i\in\B_{g'}$. Now we can replace $t$ by $tb_2^{-1}$ and $t'$ by $t'b_1$ so that now we will have ${t'}^{-1} t=r_{g'}$ and so, from (2), $M_{ i,h}=M_{fh'f^{-1}}$.
  \end{proof}

A $\Q$-valued invariant $\lambda$ of homology 3-spheres (we can assume WLOG that $\lambda (S^3)=0$) is said to be of {\em finite type} if $\lambda (\H_m )=0$ for some $m$, where $\lambda$ is extended linearly over $\H$. The {\em degree} of $\lambda$ is the smallest value of $m$ such that $\lambda (\H'_{m+1})=0$. Thus an invariant of type $m$ satisfies $\lambda (M)=\lambda (M')$ whenever $M$ and $M'$ are $A_{m+1}$-surgery equivalent. 

Conversely, according to Habiro \cite{H}, \cite{H1} if $\lambda$ is an {\em additive} invariant of homology spheres, i.e. $\lambda (M_1 \# M_2 )=\lambda (M_1 )+\lambda (M_2 )$, where $\#$ denotes connected sum, and $\lambda (M_1 )=\lambda (M_2 )$ whenever $M_1$ and $M_2$ are $A_{m+1}$-surgery equivalent, then $\lambda$ is finite-type of degree $\le m$.  This fact is proved in \cite{H} in the context of knots (Theorem 6.17), but the same arguments apply to homology spheres \cite{H1}.

For example the Casson invariant is additive and has degree  $2$.
\begin{theorem}\lbl{th.fti}
Let $\lambda$ be any rational additive invariant of finite-type. Then, for any finite $k$,  there exists (for some $g$) some $h\in\flk{k}{g}$ such that $\lambda (M_h )\not= 0$.
\end{theorem}
\begin{remark}
The same result can be proved, by essentially the same argument, for any additive finite-type invariant with values in an abelian group which has no $2$-torsion.
\end{remark}
\begin{proof}
The proof is an elaboration of the proof of Proposition \ref{th.cas}.

Let $\lambda$ be an additive invariant of finite type $d$. By the hypothesis   and Habiro's result mentioned above, there exist homology spheres $M_1$ and $M_2$ which are $A_d$-surgery equivalent such that $\lambda (M_1 )\not=\lambda (M_2 )$.  But this implies that there exists $M$ which is $A_d$-surgery equivalent to $S^3$ such that $\lambda (M)\not= 0$. We see this as follows. Let $\C$ be the semi-group, under connected sum, of oriented homology 3-spheres and let $\C/A_{d+1}$ be the set of $A_{d+1}$-surgery equivalence classes in $\C$. According to \cite{H} $\C/A_{d+1}$ is a (commutative) group under connected sum. Therefore there exists a homology sphere $M$ such that $M +M_1 =M_2$ in $\C/A_{d+1}$. Then $\lambda (M)\not= 0$, since $\lambda (M_1 )\not=\lambda (M_2 )$ and the additivity of $\lambda$ implies that  $\lambda$ induces a homomorphism $\C/A_{d+1}\to\Q$. Furthermore $M$ is $A_d$-surgery equivalent to $S^3$ because $M_1$ and $M_2$ are $A_d$-surgery equivalent and therefore, under the projection $\C/A_{d+1}\to\C/A_d$, $M=M_1 -M_2$ goes to $0$. 

By \cite{H} and Lemma \ref{lem.hab}, since $\T_g$ is a normal subgroup of $\mg$, we can write $M=M_h$ for some $h\in (\T_g )^d$ for some $g$. Now suppose, by induction, that there exists $h\in\flk{k}{g}\cap (\T_g )^d$ so that $\lambda (M_h )\not= 0$. Since $(\T_g )^d\sub\ffk{d}{g}\sub\flk{d}{g}$ we may assume $k\ge d$. Now $J_k^L (h^2 )=2J_k^L (h)$ and so, by Lemma \ref{lem.main1} and Remark \ref{rem.main1}, there is some $b\in\bar\bg\cap (\T_g )^k$ such that $J_k^L (b)=J_k^L (h^2)$. So $J_k^L (b^{-1} h^2 )=0$ and therefore $b^{-1} h^2\in\flk{k+1}{g}\cap (\T_g )^d$. But since $M_{b^{-1} h^2}=M_{h^2}$
we have $\lambda (M_{b^{-1} h^2} )=\lambda (M_{h^2} )$. We now need the following lemma.

\begin{lemma}\lbl{lem.add}
Let $\lambda$ be an invariant of homology 3-spheres of finite type $k$. Suppose $h_1\in (T_g)_{k_1}$ and  $h_2\in (\T_g )_{k_2}$, where $k_1 +k_2 >k$. Then
$$\lambda (M_{h_1 h_2})=\lambda (M_{h_1})+\lambda (M_{h_2})$$
\end{lemma}

Assuming the lemma we can conclude that $\lambda (M_{h^2})=2\lambda (M_h )\not=0$
\end{proof}
\begin{proof}[Proof of Lemma \ref{lem.add}]
Consider the linear function $\phi :\Q [\T_g ]\to\H$, mentioned above, which assigns to $h\in\T_g$ the diffeomorphism class of $M_h$. Then, according to \cite{GL2}, $\lambda\phi (I_g^{k+1})=0$. Now if $h_i\in (\T_g )^{k_i}$ then $h_i -1\in I_g^{k_i}$ and so $(h_1 -1)(h_2 -1)\in I_g^{k_1 +k_2}\sub I_g^{k+1}$. Therefore
$$\lambda(M_{h_1 h_2}-M_{h_1}-M_{h_2})=\lambda\phi (h_1 h_2 -h_1 -h_2 )=\lambda\phi ((h_1 -1)(h_2 -1))=0$$

This proves the Lemma.
\end{proof}

\section{Relating the Lagrangian filtration to the Johnson filtration}\lbl{sec.JvsL}
It is clear from the definitions that $\fk\sub\lk$ and that $\flk{\infty}{g}\sub\lk$. It seems natural to ask whether $\lk$ is generated by these two subgroups. Since $\fk$ is normal in $\mg$, this is the same as asking whether  $\lk =\fk\cd\flk{\infty}{g}$.  We will describe a recursive approach to proving (or  disproving) this conjecture, in terms of certain questions about the Johnson homomorphisms. In particular we will prove:
\begin{proposition}\lbl{prop.JvsL}
If $k=1$ or $2$ then $\lk =\fk\cd\flk{\infty}{g}$.
\end{proposition}
The recursive step  is given by:
\begin{lemma}\lbl{lem.JvsL}
Suppose that $\lk =\fk\cd\flk{\infty}{g}$. Then $\flk{k+1}{g}=\ffk{k+1}{g}\cd\flk{\infty}{g}$ if and only if
$$\im J_k\cap\ker\{\D_k (H_g )\to\D_k (H'_g )\}=J_k (\bg\cap\fk)$$
\end{lemma}
\begin{proof}
If $\lk =\fk\cd\flk{\infty}{g}$ then it is easy to see that 
$$\flk{k+1}{g}=(\fk\cap\flk{k+1}{g})\cd\flk{\infty}{g}$$
It then follows that $\flk{k+1}{g}=\ffk{k+1}{g}\cd\flk{\infty}{g}$ if and only if 
\begin{multline*}
\fk\cap\flk{k+1}{g}=\ffk{k+1}{g}\cd (\fk\cap\flk{\infty}{g})\\
=\ffk{k+1}{g}\cd (\fk\cap\bg)
\end{multline*}
since $\flk{\infty}{g}=\bg\cap\flk{1}{g}$.
 But this latter equality is equivalent to the equality
$$J_k (\fk\cap\flk{k+1}{g})=J_k (\fk\cap\bg)$$
and, finally, it is clear that 
$$J_k (\fk\cap\flk{k+1}{g})=\im J_k\cap\ker\{\D_k (H_g )\to\D_k (H'_g )\}$$
\end{proof}
For large $k$ this lemma will be difficult to apply since we do not know enough about the image of the Johnson homomorphism except when $k\le 3$. We will use the lemma for  $k\le  2$. Using recent work of Pitsch \cite{P} one may hope to settle the case $k=3$.

\begin{proof}[Proof of Proposition \ref{prop.JvsL}]
For $k=1$ this follows from the definition of $\flk{1}{g}$ and the following lemma.

\begin{lemma}
Every {\em triangular} matrix 
$\begin{pmatrix}
I & B\\0&I
\end{pmatrix}$
in $\sp (g,\Z )$ (and so $B$ is symmetric) is realized by an element of $\flk{\infty}{g}$.
\end{lemma}
The matrix is with respect to the basis $x_1,\cdots,x_g , y_1,\cdots,y_g $. $I$ denotes  the identity matrix.
\begin{proof}

\noindent (1)  If we consider the diffeomorphism obtained by a $\pm$-Dehn twist along a meridian curve representing $x_k$ we will realize the case $B=(b_{ij})$, where $b_{kk}=\pm 1$ and all other $b_{ij}=0$. 

\noi (2) If we consider the diffeomorphism which is a Dehn twist along a curve obtained by connect summing two meridians, one representing $x_k$ and the other $x_l^{\pm 1}$, ,then we realize $B$ where 
$b_{kk}=b_{ll}=-1$ and $b_{kl}=b_{lk}=\mp 1$
and all other $b_{ij}=0$. 

Note that these diffeomorphisms lie in $\bg$ since the curves along which the Dehn twists are done bound disks in $\tg$.

But now since every symmetric matrix is a linear combination of these two types, the proof is complete.
\end{proof}

For $k=2$ we will apply Lemma \ref{lem.JvsL} with $k=1$. But Lemma \ref{lem.johnson}, reinterpreted for $\bg$ instead of $\bar\bg$, interchanging $x_i$ and $y_i$, says exactly what is required by Lemma \ref{lem.JvsL} to conclude that $\flk{2}{g}=\ffk{2}{g}\cd\flk{\infty}{g}$.
\end{proof}
\ifx\undefined\bysame
	\newcommand{\bysame}{\leavevmode\hbox
to3em{\hrulefill}\,}
\fi

\end{document}